\documentclass[12pt]{amsart}
\usepackage{amssymb}
\usepackage{amsfonts}
\usepackage{latexsym}

\newcounter{cs} \stepcounter{cs} \newcounter{ds} \stepcounter{ds}

\newcommand{\casos}{\begin{itemize}}
\newcommand{\fcasos}{\end{itemize}\setcounter{cs}{1}}
\newcommand{\cas}{\item[$(\alph{cs})$]\stepcounter{cs}}

\newcommand{\subcasos}{\begin{itemize}}
\newcommand{\fsubcasos}{\end{itemize}\setcounter{ds}{1}}

\newcommand{\matriu}[1]{\left(\begin{array}{#1}}
\newcommand{\fmatriu}{\end{array}\right)} 

\newcommand{\ld}{{\rm lim}_{_{\kern-14pt\longrightarrow \kern3pt}}} 
\newcommand{\Moo}{M_{\infty}}
\newcommand{\C}{$C^*$-algebra} 
\newcommand{\Cs}{$C^*$-algebras}
\newcommand{\del}{\delta} 
\newcommand{\eps}{\epsilon}
\newcommand{\al}{\alpha}

\newcommand{\fl}{\rightarrow} 

\newcommand{\sns}{\Leftrightarrow}

\newcommand{\espai}{\vskip 0.4cm}

\newtheorem{lem}{Lemma}[section] 
\newtheorem{corol}[lem]{Corollary}
\newtheorem{theor}[lem]{Theorem} 
\newtheorem{prop}[lem]{Proposition}
\newtheorem{rema}[lem]{Remark}

\newcommand{\End}{{\rm End}}

\begin{document}
\title[Lifting units]{Lifting units modulo exchange ideals and \boldmath 
$C^*$-algebras with real rank zero}
\author{Francesc Perera}
\thanks{Partially supported by the TMR EU-Programme, by the DGESIC (Spain) 
  and by the Comissionat per Universitats i Recerca de la Generalitat
  de Catalunya.}
\address{Department of Mathematics, University of Copenhagen,
  Universitetsparken 5, 2100 Copenhagen \O, Denmark}
\email{perera@math.ku.dk}
\address{Departament de Matem\`atiques, Universitat Aut\`onoma de
  Barcelona, 08193, Bellaterra (Barcelona), Spain} 
\email{perera@mat.uab.es} 
\date{} \dedicatory{} \commby{}
\keywords{Exchange ring, $C^*$-algebra, real rank, refinement,
  separativity, $K$-Theory, index.}

\begin{abstract}
Given a unital ring $R$ and a two-sided ideal $I$ of $R$, we consider
the question of determining when a unit of $R/I$ can be lifted to a
unit of $R$. For the wide class of separative exchange ideals $I$, we
show that the only obstruction to lifting invertibles relies on a
$K$-theoretic condition on $I$. This allows to extend previously known 
index theories to this context. Using this we can draw consequences 
for von Neumann regular rings and \Cs\ with real rank zero.
\end{abstract}

\maketitle

\section*{Introduction}

The problem of lifting units from a quotient of a ring $R$ modulo a
two-sided ideal $I$ has been of interest in several instances. The
first extension of the classical index theory for Fredholm operators
on a Hilbert space was directed to von Neumann algebras (see
\cite{br1}, \cite{br2} and \cite{ol}). The class of self-injective
rings and that of Rickart \Cs\ satisfying certain comparability
conditions were considered by Menal and Moncasi (see \cite{memo}). The
general case for Rickart \Cs\ was studied by Ara in \cite{pkt}. In all 
of the above cases, the extent to which a unit from a quotient $R/I$ can
be lifted to a unit of $R$ is measured by a condition of $K$-theoretic
nature, namely the vanishing of the connecting index map in $K$-Theory.

Our aim here is to consider the class of exchange ideals of unital
rings, which is known to contain both (not necessarily unital) von
Neumann regular rings and \Cs\ with real rank zero. In fact, the
exchange \Cs\ are exactly those having real rank zero (see
Theorem 7.2 in \cite{agop}). A second unifying principle on which we
will rely is that of separative cancellation of finitely generated
projective modules, which can be regarded as a weak cancellation
property. Separative unital exchange rings provide a framework in
which a number of outstanding open problems are known to have solutions
(see \cite{agopl}, \cite{agop}). Moreover, this weak cancellation
condition holds widely (for instance, for all the known classes of
regular rings - see \cite{agop}, and also for the known classes of extremally
rich \Cs\ - see \cite{bpp}); it is therefore regarded as a condition
that might hold for all exchange rings.

Our main objective is to prove that if $I$ is a separative exchange
ideal of a unital ring $R$, then the index is the only obstruction to
lifting units modulo $I$, thus providing a common setting in which the 
results mentioned above can be handled. Hence we derive some
consequences for both regular rings and \Cs\ with real rank zero (and
their multiplier algebras). In
order to develop our index theory, we benefit from results and
techniques from \cite{agor}, where a detailed analysis of elementary
transformations on invertible matrices over {\it unital} exchange
rings was carried out. Our present context is, however, different in
that we deal with invertible matrices over unital rings, which are
diagonal modulo an exchange ideal. It is  remarkable that, in a different
direction, it has been established by Brown and Pedersen
(\cite{bpp}, Theorem 5.2) that the index is also the only obstruction
to lifting invertibles modulo separative, extremally rich ideals of
unital \Cs.

We now fix some notations. As a general rule, $R$ will stand 
for a unital ring, whereas we shall use $I$ to denote a nonunital
ring, generally sitting inside $R$ as a two-sided ideal. An {\it
  elementary matrix} is a matrix of the form {\boldmath $1$}$+re_{ij}$,
where {\boldmath $1$} is an identity matrix, $e_{ij}$ is one of the
usual matrix units (with $i\neq j$), and $r\in R$. We will denote by $E_n(R)$ the
subgroup of $GL_n(R)$ generated by the elementary matrices. If $x,y\in
M_n(R)$, then we use $x\oplus y$ to denote the matrix
$\left(\begin{array}{cc} x & 0 \\ 0 &  y \end{array}\right)$, and we
will denote by $1_n$ the unit of $M_n(R)$.

\section{Preliminary results}

Let $M$ be a right $R$-module. We say that $M$ satisfies the {\it
finite exchange property} (see \cite{cwj}) if for every right
$R$-module $A$ and any decompositions 
$$A=M'\oplus N=A_1\oplus\ldots\oplus A_n$$
with $M'\cong M$, there exist submodules $A_i'\subseteq A_i$ (which
are in fact direct summands, by the modular law) such that
$$A=M'\oplus A_1'\oplus \ldots \oplus A_n'.$$

Following \cite{war}, we say that $R$ is an {\it exchange ring} provided
that $R_R$ satisfies the finite exchange property. This notion is
right-left symmetric (see \cite{war}, Corollary 2). In \cite{nic},
Theorem 2.1, it is proved that $M$ has the finite exchange property if
and only if $\End (M)$ is an exchange ring. Also, a useful ring-theoretic
characterization was provided independently by Goodearl and Nicholson:
\begin{lem}
\label{pr1}
{\rm (\cite{gw}, Theorem 2.1 in \cite{nic})}
A unital ring $R$ is an exchange ring if and only if for every element 
$a\in R$ there exists an idempotent $e\in aR$ such that $1-e\in (1-a)R$.\qed
\end{lem}

This characterization motivated the notion of an exchange ring for rings
without unit (see \cite{pext}). Namely, a (possibly non-unital)
ring $I$ is said to be an {\it exchange ring} if for each $x\in I$,
there exist an idempotent $e\in I$ and elements $r,s\in I$ such that
$e=xr=x+s-xs$. As proved in Lemma 1.1 in \cite{pext}, the ring $I$ is
exchange if and only if, whenever $x\in I$ and $R$ is a unital ring
containing $I$ as a two-sided ideal, then there exists an idempotent
$e\in xI$ such that $1-e\in (1-x)R$. Of course, the notions of unital
exchange and non-unital exchange agree if the ring $I$ has a unit.

Since we will usually have a non-unital exchange ring $I$ which is an
ideal of a unital ring $R$, we will adopt the terminology in
\cite{pext} and say that, in this context, ``$I$ is an exchange ideal
of $R$''.

The class of exchange rings is pleasantly large: it includes regular
rings, $\pi$-regular rings, semiperfect rings (which are exactly the
semilocal exchange rings), right self-injective rings (see
Remark 2.9 (b) in \cite{agor}) and \Cs\ with real rank zero
(by Theorem 7.2 in \cite{agop}). It is closed under natural ring
constructions. For example, if $I$ is an exchange ideal of $R$ and
$e\in R$ is an idempotent, then $eIe$ is an exchange ideal of $eRe$
(\cite{pext}, Proposition 1.3 and also \cite{nic}, Proposition
1.10). Also, if $I$ is an exchange ring, then $M_n(I)$ is an
exchange ring for all $n\geq 1$ (by Theorem 1.4 in \cite{pext}). The
behaviour of exchange rings under extensions is characterized by an
idempotent-lifting condition (\cite{pext}, Theorem 2.2): if $I$ is an
ideal of a (possibly non-unital) ring $L$, then $L$ is exchange if and
only if $I$ and $L/I$ are exchange and idempotents can be lifted
modulo $I$.

Let $I$ be an ideal of a unital ring $R$. We denote by $FP(I,R)$ the
class of all finitely generated projective right $R$-modules $P$ such that
$P=PI$, and we define $V(I)$ to be the set of all isomorphism classes
of elements from $FP(I,R)$. Note that $V(I)$ becomes an abelian monoid
under the operation $[P]+[Q]=[P\oplus Q]$. Even though $V(I)$ involves 
the unital ring $R$, it can be shown that it only depends on the ring
structure of $I$ (see, for example, \cite{memo}, or also
\cite{ros}). It will be
sometimes convenient to use an alternate description of $V(I)$ via idempotents
(see \cite{ros}), so we identify $V(I)$ with the set of equivalence
classes of idempotents in $\Moo(I)$, the non-unital ring of
$\omega\times\omega$ matrices with only finitely many nonzero entries
from $I$. We will use $[e]$ to 
indicate the class of $e$ in $V(I)$. Viewing $I$ inside $R$, we can
also identify $V(I)$ with $\{[e]\in V(R)\mid e \mbox{ an idempotent in 
  } \Moo(I)\}$.

Let $M$ be an (abelian) monoid. We can order $M$ by using the
so-called algebraic ordering: for $x,y\in M$, write $x\leq y$ provided
that there exists $z\in M$ such that $x+z=y$. If $S$ is a submonoid of 
$M$ with the property that if $x\leq y$ and $y\in S$, then $x\in S$,
then $S$ is called an {\it order-ideal} of $M$. For example, if $I$ is 
a two-sided ideal of a ring $R$, then $V(I)$ is an order-ideal of
$V(R)$. Notice that if $e$ and $f$ are idempotents in $M_n(R)$ for
some $n$, then $[e]\leq [f]$ is equivalent to saying that $eM_n(R)$ is
isomorphic to a direct summand of $fM_n(R)$.

Let $M$ be a monoid, and let $S$ be an order-ideal of $M$. We say that
$M$ has {\it refinement with respect to $S$} if whenever
$x_1+x_2=y_1+y_2$ in $M$ with at least one of $x_1$, $x_2$, $y_1$,
$y_2$ in $S$, then there exist elements $z_{ij}\in M$ such that
$x_i=z_{i1}+z_{i2}$ and $y_i=z_{1i}+z_{2i}$, for $i=1,2$.

Observe that if $S=M$ then we are recovering the usual definition of a 
refinement monoid (see, for example \cite{agop}). Also, we see from the
definition that in particular, $S$ is a refinement monoid. If $I$ is
an exchange ring, it is known that $V(I)$ is a refinement monoid (see
\cite{pext}, Proposition 1.5, and also \cite{agop}, Proposition
1.2). If now $I$ is an exchange ideal of a unital ring $R$,
still some refinement persists in $V(R)$, as the following lemma
shows:
\begin{lem}
\label{pr2}
Let $I$ be an exchange ideal of a unital ring $R$. Then $V(R)$ has
refinement with respect to $V(I)$.
\end{lem}
{\em Proof:} Let $[A]+[B]=[C]+[D]$, with $[A]$, $[B]$, $[C]$, $[D]\in
V(R)$, and assume, for example, that $[A]\in V(I)$, that is, $A\in
FP(I,R)$. Then $\End (A)$ is an exchange ring, so that $A$ has the
finite exchange property, and hence we may use the proof of
Proposition 1.2 in \cite{agop}.\qed
\espai
We say that an abelian monoid $M$ is {\it separative} provided that
whenever $a+a=a+b=b+b$ in $M$, then $a=b$. Equivalently, $M$ is
separative if the following weak cancellation condition holds: if
$a+c=b+c$ and $c\leq na, nb$ for some $n$, then $a=b$ (see
Lemma 2.1 in \cite{agop}). Accordingly, we call a ring $R$ {\it separative} if
$V(R)$ is a separative monoid (see \cite{agop}). The following lemma
was stated in \cite{agop}, Lemma 4.4, for full refinement monoids, and the
proof used there can be used in our present context almost
entirely. Since it will be an essential result later, we just indicate
the major steps that lead to the conclusion.

\begin{lem}
\label{pr3}
Let $M$ be a monoid and let $S$ be a separative order-ideal such that $M$ has
refinement with respect to $S$. If $a+e=b+e$ for $a,b\in M$ and $e\in
S$, and $e\leq na,nb$ for some $n$, then $a=b$.
\end{lem}
{\em Proof:} Since $M$ has refinement with respect to $S$, and $e\in S$ 
with $e\leq
na$, we can decompose $e=\sum\limits_{i=1}^n e_i$, with $e_i\leq a$
for all $i$. Hence we may assume that $e\leq a$, and similarly $e\leq
b$.

By refinement with respect to $S$ we get decompositions $a=a_1+a_2$,
$b=a_1+b_2$, and $e=b_2+c_2=a_2+c_2$. Note that $c_2\in S$ and that
$c_2\leq e\leq a=a_1+a_2$, hence we may use refinement with respect to 
$S$ again. It is not difficult to see that we can arrange the
resulting decompositions and change notation in such a way that
$c_2\leq a_2,b_2$. Now separativity in $S$ entails $a_2=b_2$, so $a=b$.\qed
\espai
Recall that an element $x$ in a ring $R$ is called {\it von Neumann
  regular} if there exists $y\in R$ such that $x=xyx$. If $y$ can be
chosen to be a unit, then $x$ is called {\it unit-regular}. In this
case $x=(xy)y^{-1}$ is a product of an idempotent with a unit.

\begin{prop}
\label{pr4}
Let $I$ be a separative exchange ideal of a unital ring $R$, and let
$d\in I$ be such that $dR=(1-p)R$ and $Rd=R(1-q)$, for some
idempotents $1-p,1-q\in I$. If $RpR=RqR=R$ then $[p]=[q]$ in
$V(R)$. In particular, $d$ is unit-regular.
\end{prop}
{\em Proof:}
Note that $V(I)$ is a separative monoid and that $V(R)$ has refinement
with respect to $V(I)$ by Lemma \ref{pr2}. From the outset we 
evidently have that $[1-p]=[1-q]$ in $V(I)\subseteq  V(R)$. Also,
since $RpR=RqR=R$, we get that $[1-p]\leq n[p], n[q]$ for some $n\in
\mathbb N$. Now, in $V(R)$, we have that
$[p]+[1-p]=[q]+[1-p]$. By Lemma \ref{pr3}, we get that $[p]=[q]$ in
$V(R)$, so that $pR\cong qR$. Now, note that $R/dR=R/(1-p)R\cong
pR\cong qR={\rm r.ann }(d)$, whence $d$ is unit-regular, by the proof
of Theorem 4.1 in \cite{vnrr}.\qed

\section{Index theory for exchange rings}

\begin{lem}
\label{lif1}
{\rm (}cf. \cite{agor}{\rm, Lemma 2.1)}
Let $I$ be an exchange ideal of a unital ring $R$, and let
$e_1,e_2\in R$ be idempotents such that $e_1\in I$. Then, there exists
an idempotent $e\in e_1R+e_2R$ such that $[e_i]\leq [e]$ in $V(R)$,
for all $i$. In particular, $ReR=Re_1R+Re_2R$.
\end{lem}
{\em Proof:}
Considering that $\End(e_1R)=e_1Re_1=e_1Ie_1$ is a unital exchange
ring, we have that $e_1R$ has the finite exchange property, and hence
we get decompositions $e_2R=A\oplus B$ and $(1-e_2)R=A'\oplus B'$ such
that $R=e_1R\oplus A\oplus A'$. Then, choose $e\in R$ such that
$eR=e_1R\oplus A$, and proceed as in the proof of Lemma 2.1 in
\cite{agor}.\qed
\espai
In the next technical lemmas, we will be involved
with performing several elementary row and column operations on an
invertible $ 2\times 2$ matrix over a ring $R$. These operations will
follow the lines of 2.3-2.7 in \cite{agor}. However, our ring $R$ won't be
exchange, so we cannot apply the results in \cite{agor} directly, and
hence some different procedure is needed.

Let $I$ be a two-sided ideal in a unital ring. We shall denote by
$\pi:R\fl R/I$ the natural quotient map. For any $n>1$, let
$E_n(I)$ be the subgroup of $E_n(R)$ generated by the elementary matrices
$1_n+re_{ij}$ for $r\in I$ and $i\neq j$. Note that $\pi (\eps)=1_n$
for all $\eps\in E_n(I)$. Observe also that multiplying a matrix
$\al\in M_n(R)$ on the left or right by any matrices from $E_n(I)$ does 
not change $\pi (\al)$.
\begin{lem}
\label{lif4a}
Let $I$ be an exchange ideal of a unital ring $R$. Let
$\al=\left(\begin{array}{cc} a & b \\ c & d \end{array}\right)\in
GL_2(R)$ such that $b,c\in I$.
\casos
\cas
There exist $\beta\in E_2(I)$ and an idempotent $1-h\in I$ such that
$\al\beta =\left(\begin{array}{cc} a' & b' \\ c' & d'
  \end{array}\right)$, with $c'\in Rc$, $c'R=(1-h)R$, $d'R=hR$ and $RhR=R$.
\cas
There exist $\gamma\in E_2(I)$ and an idempotent $1-k\in I$ such that
$\gamma\al=\left(\begin{array}{cc} a'' & b'' \\ c'' & d''
  \end{array}\right)$, with $b''\in bR$, $Rb''=R(1-k)$, $Rd''=Rk$ and $RkR=R$.
\fcasos
\end{lem}
{\em Proof:}
$(a)$. First note that the row $(c,d)$ is right unimodular, so
$cR+dR=R$. Hence there exist $x,y\in R$ such that $cx+dy=1$. Now,
since $c\in I$ and $I$ is exchange, there exists an idempotent $e\in
cxR\subseteq cR$ such that $1-e\in dR$. Write $e=cr$, with $re=r$ and
$1-e=ds$, with $s(1-e)=s$. As in \cite{agor}, Lemma 2.4, we multiply
$\al$ on the right by $\al_1\al_2$, where
$\al_i\in E_2(I)$ are defined as follows:
$$\al_1=\left(\begin{array}{cc} 1 & 0
    \\ -sc & 1 \end{array}\right),\quad\al_2=\left(\begin{array}{cc} 1 & -rd
    \\ 0 & 1 \end{array}\right).$$
We thus obtain as last row: $(ec,(1-e)d)$. As in the proof of Lemma
2.4 in \cite{agor} we get $ec\in cRc$, $(1-e)d\in dRd$, and $R=ecR\oplus
(1-e)dR$. Set $w=ec+(1-e)d$, and note that $ewr+(1-e)ws=1$, hence
$RwR=R$. Similar to the proof of Corollary 2.5 in \cite{agor}, we want to
put $w$ in the $(2,2)$ position. By the exchange property again
(applied to $ewr\in I$), we get an idempotent $f\in ewrR\subseteq ewR$
such that $1-f\in (1-ewr)R\subseteq (1-e)wR$. Write $f=eww_1$, with
$w_1f=w_1$, and $1-f=(1-e)ww_2$, with $w_2(1-f)=w_2$, and set
$f_1=ww_1e\in wR$, and $f_2=ww_2(1-e)\in wR$. Note that $f_i$ are
idempotents with $f_1\in I$, and that $f_1R\cong fR$, whereas
$f_2R\cong (1-f)R$. By Lemma \ref{lif1}, there is an idempotent $g\in
f_1R+f_2R\subseteq wR$ such that $RgR=Rf_1R+Rf_2R$. Since $Rf_1R=RfR$
and $Rf_2R=R(1-f)R$, we see that $RgR=R$.

Write $g=ww'$, for some $w'\in R$, and multiply $\al\al_1\al_2$ on the
right by $\beta_1\beta_2$, where $\beta_i\in E_2(I)$ are defined as:
$$
\beta_1=\left(\begin{array}{cc}
1 & rc \\ 0 & 1
\end{array}\right),\quad
\beta_2=\left(\begin{array}{cc}
1 & 0 \\ -w'ec & 1
\end{array}\right).
$$
This gives as last row: $((1-g)ec,w)$. At this point we start the 
first part of the procedure again with the current last row, so that after
right multiplication by two matrices $\gamma_i\in E_2(I)$ for $i=1,2$,
we get as last row $(c',d')$, with $c'\in 
(1-g)ecR(1-g)ec\subseteq (1-g)Rc$, and $R=c'R\oplus d'R$. According to
the direct sum decomposition, we see that there exists an idempotent
$1-h\in I$ such that $c'R=(1-h)R$ and $d'R=hR$. Since $g(1-h)R=gc'R=0$ 
and $RgR=R$, we conclude that $RhR=R$. Finally, set
$\beta=\al_1\al_2\beta_1\beta_2\gamma_1\gamma_2$.

(b). Since $b\in I$ and $I$ is exchange, we can perform the transpose version 
of the process carried out in $(a)$ to the last column. Hence we
obtain matrices $\al_i'$, $\beta_i'$ and $\gamma_i'$ in $E_2(I)$, for
$i=1,2$, such that after left multiplication by
$\gamma:=\gamma_2'\gamma_1'\beta_2'\beta_1'\al_2'\al_1'$ we get a
matrix $\left(\begin{array}{cc} a'' & b''\\ c'' & d''
  \end{array}\right)$ and an idempotent $1-k\in I$ satisfying the
desired properties.\qed
\begin{lem}
\label{lif4}
Let $I$ be a separative exchange ideal of a unital ring $R$, and let
$\al=\left(\begin{array}{cc} a & b \\ c & d \end{array}\right)\in
GL_2(R)$ such that $b,c\in I$ and $d-1\in I$. Then there exist units
$a',u\in GL_1(R)$, and matrices $\beta$, $\gamma$ and $\eps$ in
$E_2(R)$ such that
$$\gamma\al\beta(1\oplus u^{-1})\eps=a'\oplus 1,$$ and
$\pi(a')=\pi (au^{-1})$. In particular, $[\al]=[a'u]$ in $K_1(R)$.
\end{lem}
{\em Proof:}
We remark again that $\pi(\al)$ remains unchanged after right or left
multiplication by matrices from $E_2(I)$. Thus we apply Lemma
\ref{lif4a} $(a)$, and without loss of generality there is an
idempotent $1-h\in I$ such that $cR=(1-h)R$, $dR=hR$ and $RhR=R$.

Now we proceed as in the proof of Lemma 2.7 in \cite{agor}, so that we
move $c$ to the $(1,2)$ position. This is achieved after right and
left multiplication by the signed permutation matrix
$\sigma=\left(\begin{array}{cc} 
    0 & 1 \\ -1 & 0 \end{array}\right)$, which is a product of three
elementary matrices. Hence, we get:
$$\al':=\sigma\al\sigma=\left(\begin{array}{cc}
-d & c \\ b & -a
\end{array}\right).$$

By Lemma \ref{lif4a} $(b)$, there exist $\gamma'\in E_2(I)$ and an
idempotent $1-q\in I$ such that
$\gamma'\al'= \left(
\begin{array}{cc}
a' & b' \\ c' & d'
\end{array}\right)$, with $b'\in cR=(1-h)R$, $Rb'=R(1-q)$ and
$RqR=R$. Note that $\pi(a')=-1$. Since $b'$ is regular, there exists
an idempotent $1-p\in I$ such that $(1-p)R=b'R$. Using that
$h(1-p)R=hb'R\subseteq hcR=0$ and $RhR=R$, we conclude that
$RpR=R$. Hence, by Proposition \ref{pr4} (and since $I$ is separative,
by hypothesis), $b'$ is unit-regular. Write $b'=fu$, where $f^2=f\in
I$ and $u\in GL_1(R)$.

Now move the element $b'$ to the $(2,2)$ position, multiplying
$\gamma'\al'$ on the left by $\sigma^{-1}$. Then we multiply
$\sigma^{-1}\gamma'\al'$ on the right by $\lambda
=\left(\begin{array}{cc} 1 & 0 \\ 0 & u^{-1}\end{array}\right)$. The
matrix we obtain is $\al'':=\left(\begin{array}{cc} r & s \\ t & f
  \end{array}\right)$, where $r\in I$, and $\pi(t)=-1$.

We want to use now Lemma 2.3 in \cite{agor}, in order to get $1$ in the
$(2,2)$ position and zeros elsewhere in the last row and
column. First, multiply $\al''$ on the right by 
$\eps_1=\left(\begin{array}{cc} 1 & 0 \\ -ft & 1 \end{array}\right)\in 
E_2(I)$, 
so we get as last row $((1-f)t, f)$. Since
this row is right unimodular, we have that $(1-f)tR+fR=R$, and hence
$(1-f)tR=(1-f)R$. Choose $v\in R(1-f)$ such that 
$(1-f)tv=(1-f)$. Since $\pi(t)=-1$, we have that $\pi(v)=-1$. After
right multiplication by $\eps_2=\left(\begin{array}{cc} 1 & v \\ 0 & 1
  \end{array}\right)$, our second row becomes $((1-f)t,1)$. Denote by
$z'$ the element in the $(1,2)$ position of the resulting matrix.

Finally, set $\mu_1=\left(
  \begin{array}{cc}
1 & -z' \\ 0 & 1
  \end{array}
\right)$ and $\mu_2=\left(
  \begin{array}{cc}
1 & 0 \\ -(1-f)t & 1
  \end{array}
\right)$, the last routine matrices. Observe that 

$$\pi(\mu_1)=\left( 
  \begin{array}{cc}
1 & \pi(-z') \\ 0 & 1
  \end{array}
\right)
\mbox{ and that }\pi(\mu_2)=\left(\begin{array}{cc} 1 & 0 \\ 1 &
    1\end{array}\right).$$ After multiplying $\al''\eps_1\eps_2$ on
the left by $\mu_1$ and on the right by $\mu_2$, we get a matrix of the form
$a'\oplus 1$, where $a'\in R$. We need to compute $\pi(a')$. By all
the calculations performed so far, we have that:

$$\pi(a')\oplus 1=\pi
(\mu_1\sigma^{-1}\gamma'\sigma\al\sigma\lambda\eps_1\eps_2\mu_2)=$$
$$\left(\begin{array}{cc} 1 & \pi(-z')\\ 0 &
    1\end{array}\right) 
\left(\begin{array}{cc} \pi(a) & 0 \\ 0 & 1 \end{array}\right)
\left(\begin{array}{cc} 0 & 1 \\ -1 & 0 \end{array}\right)
\left(\begin{array}{cc} 1 & 0 \\ 0 & \pi(u^{-1}) \end{array}\right)
\left(\begin{array}{cc} 1 & -1 \\ 0 & 1 \end{array}\right)
\left(\begin{array}{cc} 1 & 0 \\ 1 & 1 \end{array}\right).$$

After computing the right-hand side of the above equality, we obtain:
 $$\pi(a')\oplus 1=\left(  
  \begin{array}{cc}
\pi(au^{-1}) & \pi(au^{-1})-\pi(z') \\ 0 & 1
  \end{array}
\right),$$ whence we get that $\pi(a')=\pi(au^{-1})$, as
desired.\qed
\espai
Let $R$ be a ring, and let $I$ be a two-sided ideal. We define the
{\it Fredholm elements relative to $I$} as the set
$F(I,R)=\pi^{-1}(GL_1(R/I))$.

Note that $F(I,R)$ is a multiplicative subsemigroup of $R$, such that
$GL_1(R)+I\subseteq F(I,R)$. Observe that if $R$ has stable rank one, then
$GL_1(R)+I=F(I,R)$. Denote by $\del:K_1(R/I)\fl K_0(I)$ the
connecting map in algebraic $K$-Theory, and recall that there is an
exact sequence
$$K_1(R)\stackrel{\pi_1}{\fl} K_1(R/I)\stackrel{\del}{\fl} K_0(I)\fl
K_0(R)\fl K_0(R/I)$$ 
(\cite{ros}, Theorem 2.5.4). As usual, cf. \cite{bpi}, 6.1, we define the {\it
  index map} as the semigroup homomorphism $${\rm index}: F(I,R)\fl
K_0(I),$$ given by the rule ${\rm index}(x)=\del([\pi(x)])$.

We are now in position to prove our main result:
\begin{theor}
\label{lif5}
Let $I$ be a separative exchange ideal of a unital ring $R$. Let $x\in
R$ be a Fredholm element relative to $I$. Then there exists
$y\in GL_1(R)$ such that $x-y\in I$ if and only if ${\rm index}(x)=0$. In
this case, for any $\al\in K_1(R)$ that is mapped to $[\pi(x)]$, we may find
$y\in GL_1(R)$ such that $[y]=\al$ and $x-y\in I$.
\end{theor}
{\em Proof:}
Assume that $\pi(x)$ can be lifted to a unit of $R$. Then there
exists $y\in GL_1(R)$ such that $x=y+b$, with $b\in I$. Hence $${\rm
  index}(x)=\del ([\pi(x)])=\del ([\pi(y)])=\del \pi_1([y])=0,$$
by exactness.

Conversely, suppose that ${\rm index}(x)=0$. Again by exactness this
means that there exists $k\in\mathbb N$ and $y_1\in GL_k(R)$ such that
$[\pi(y_1)]=[\pi(x)]$. Hence, if $m=2^n$ is large enough, there is
$\pi(z)\in E_m(R/I)$ such that $\pi (x)\oplus 1_{m-1}=\pi
(z)(\pi(y_1)\oplus 1_{m-k})$. We may clearly assume that in fact $z\in 
E_m(R)$. Set $w_1:=z(y_1\oplus 1_{m-k})$, and note that $[w_1]=[y_1]$
in $K_1(R)$. Denote by $(a_{ij})$ the
entries of $w_1$, and observe that $a_{ij}\in I$ whenever $i\neq j$,
that $a_{ii}-1\in I$ for all $i\neq 1$, and that $a_{11}-x\in
I$.

We apply Lemma \ref{lif4}, replacing $I$ and $R$ by $M_{m/2}(I)$ and
$M_{m/2}(R)$. (Notice that $M_{m/2}(I)$ is a separative exchange ideal of
$M_{m/2}(R)$.) Thus we obtain matrices $u,w_2'\in GL_{m/2}(R)$ such
that $\pi(w_2'u)=(\pi(x)\oplus 1_{m/2-1})$, and $[w_1]=[w_2'u]$ in
$K_1(R)$. Let $w_2:=w_2'u$. A recursive procedure shows that we get
$w_n\in GL_1(R)$ such that $\pi(w_n)=\pi(x)$, and $[w_n]=[y_1]$ in
$K_1(R)$.\qed
\begin{rema}
\label{manola}
{\rm Notice that, according to Theorem 3.5 in \cite{ags}, every ring has a 
largest exchange ideal with respect to the inclusion (which might be
the zero ideal).}
\end{rema}
\begin{corol}
Let $R$ be a (unital) separative exchange ring such that $K_0(R)=0$,
and let $I$ be an ideal of $R$. Then the units of $R/I$ lift to those
of $R$ if and only if $K_0(I)=0$.
\end{corol}
{\em Proof:}
Clearly, if $K_0(I)=0$, then the connecting map $\del$ vanishes and
hence Theorem \ref{lif5} applies. Conversely, suppose that the units of 
$R/I$ lift to units of $R$. By \cite{agor}, Theorem 2.8, the natural
map $GL_1(R/I)\fl K_1(R/I)$ is surjective. It follows then that the
map $\pi_1:K_1(R)\fl K_1(R/I)$ in $K$-Theory is surjective. Therefore
we get again that $\del=0$, and since $K_0(R)=0$, we conclude by
exactness that $K_0(I)=0$.\qed
\espai
The previous corollary applies to the case where $R$ is a purely
infinite, right self-injective ring and recovers some results of Menal 
and Moncasi (\cite{memo}). If $R$ is a right self-injective
ring, then by Theorem 1.22 in \cite{vnrr}, the quotient $R/J(R)$ is
regular and right self-injective, where $J(R)$ is the Jacobson radical
of $R$. We say that a right self-injective ring is {\it purely
  infinite} if $R/J(R)$ is a purely infinite regular ring. (According
to \cite{vnrr}, a regular ring $R$ is {\it purely infinite} if $R\cong
R\oplus R$.) These rings are known to be
exchange and separative. Further, $K_0(R)=0$, by the proof of
Corollary 3.6 in \cite{pext}.

\section{\Cs\ with real rank zero}

Let $A$ be a \C. In order to distinguish between the algebraic and the 
topological $K_1$-groups of $A$, and according to more common usage, we
will denote by $K_1^{{\rm alg}}(A)$ the algebraic $K_1$-group of $A$,
and we shall use $K_1(A)$ to denote the topological $K_1$-group (see
\cite{bl}, Definition 8.1.1, \cite{wo}, Definition 7.1.1). It is
known that there is a natural surjective homomorphism $\gamma:K_1^{{\rm
    alg}}(A)\fl K_1(A)$ (see, for example, \cite{agor}). Since
idempotents in $M_n(A)$ are equivalent to projections (e.g., \cite{bl}),
we may identify $V(A)$ with the abelian monoid of Murray-von Neumann
equivalence classes of projections arising from $\Moo(A)$.

Recall that a (unital) \C\ $A$ has {\it real rank zero} provided that
every self-adjoint element can be approximated arbitrarily well by
self-adjoint, invertible elements. Other characterizations, including 
the original definition, may be found in \cite{bpf}. If $A$ is
non-unital, then $A$ has real rank zero if and only if the minimal
unitization $\widetilde{A}$ of $A$ (see \cite{wo}) has real rank
zero. Since the \Cs\ that are exchange rings are precisely those
having real rank zero, we see that if $A$ is a \C\ with real rank
zero, then $V(A)$ is a refinement monoid (see also \cite{zhr}, Theorem 
5.3).

Brown and Pedersen have introduced in \cite{bpp} the concept of {\it weak
cancellation} for \Cs, meaning that if $p$ and $q$ are projections in $A$
that generate the same closed ideal $I$ of $A$, and $[p]=[q]$ in
$K_0(I)$, then they are (Murray-von Neumann) equivalent in $A$. If
this property holds for $M_n(A)$, for all $n$, then $A$ has {\it
  stable weak cancellation}. Notice that $A$ has stable weak
cancellation if and only if $A$ is separative. In fact, if $A$ has
real rank zero, then $A$ has weak cancellation if and only if $A$ is
separative (that is, the property of weak cancellation is
stable). This follows using Proposition 2.8 in \cite{agop} and the fact
that $V(A)$ is a refinement monoid.

The property of (stable) weak cancellation is shown to hold widely
within the class of extremally rich \Cs\ (see \cite{bpe}),
including those that have real rank zero (see \cite{bpp}, Theorem
2.11). As for the case of exchange rings, there are no examples
known of extremally rich \Cs\ without weak cancellation (\cite{bpp},
Remark 2.12).

Let $A$ be a \C, and let $I$ be a closed, two-sided ideal of
$A$. Denote by $\partial:K_1(A/I) \fl K_0(I)$ the connecting map in
topological $K$-Theory (see, e.g., \cite{bl}, Definition 8.3.1,
\cite{wo}, Definition 8.1.1). We then define the index of a Fredholm
element $x$ (relative to $I$) as ${\rm index}(x)=\partial([\pi(x)])$,
where $\pi:A\fl A/I$ is the natural projection map. Now, since we have
that $\partial\gamma=\del$, where $\del:K_1^{{\rm alg}}(A/I)\fl
K_0(I)$ is the algebraic connecting map, we see that the two possible
definitions of algebraic and topological indices for Fredholm elements
coincide. From the observations made, it is clear that we can apply
the result in the previous section to get the following theorem (which
has been independently obtained by L.G. Brown [unpublished]). We
remark that for extremally rich ideals with weak cancellation the same
conclusion holds, as shown in \cite{bpp}, Theorem 5.2.
\begin{theor}
\label{c1}
Let $A$ be a \C\ and let $I$ be a closed ideal of $A$ with real rank
zero and weak cancellation. Let $x\in A$ be a Fredholm element
relative to $I$. Then there exists $y\in
GL_1(A)$ such that $x-y\in I$ if and only if ${\rm index}(x)=0$. In
this case, for any $\al\in K_1(A)$ that is mapped to $[\pi(x)]$, we may find
$y\in GL_1(A)$ such that $[y]=\al$ and $x-y\in I$.
\end{theor}
{\em Proof:}
As in the proof of Theorem \ref{lif5}, if there exists $y\in GL_1(A)$ such
that $x-y\in I$, then ${\rm index}(x)=0$.

Conversely, if ${\rm index}(x)=\partial([\pi(x)])=0$, then since
$\gamma([\pi(x)])=[\pi(x)]$ we have in fact that $\delta([\pi(x)])=0$, 
and so Theorem \ref{lif5} applies.\qed
\begin{rema}
{\rm As for the case of exchange rings (see Remark \ref{manola}), any
  \C\ has a largest closed ideal of real rank zero, a description of which is 
  given in \cite{bpid}, Theorem 2.3.}
\end{rema}
We now give some applications to the multiplier algebras $M(A)$ of
\Cs\ $A$ with real rank zero. Multiplier algebras of \Cs\ are relevant 
objects (since they can be used, for instance, to parametrize
extensions) that have been intensively studied in the last years (to
cite a few examples, among many others, see \cite{apt}, \cite{ell},
\cite{brcan}, \cite{linpams}, \cite{zh}, \cite{gkth}, \cite{per}). For
the basic facts concerning multipliers see, for example, Chapter 2 in
\cite{wo}.

The proof of the following corollary is derived entirely as the proof
of Corollary 5.8 in \cite{bpp}, using Theorem \ref{c1} instead of
\cite{bpp}, Theorem 5.2. For any unital \C\ $A$, we use $U(A)$ to
denote the unitary group of $A$, whereas $U_0(A)$ stands for the
connected component of the identity in $U(A)$.
\begin{corol}
\label{c2}
Let $A$ be a $\sigma$-unital and stable \C. Suppose that $A$ has real
rank zero and weak cancellation. Then there is a short exact sequence
of groups
$$0\fl U_0(M(A)/A)\fl U(M(A)/A)\fl K_0(A)\fl 0.\qed$$
\end{corol}
\begin{corol}
\label{c3}
Let $A$ be a $\sigma$-unital and stable \C, with real rank zero and
weak cancellation. Then the following are equivalent:
\casos
\cas
The units of $M(A)/A$ lift to units of $M(A)$;
\cas
$K_0(A)=0$;
\cas
$U(M(A)/A)$ is connected.
\fcasos
\end{corol}
{\em Proof:}
Since $A$ is $\sigma$-unital and stable, we have that
$U(M(A))$ is connected (see \cite{mi}, \cite{ch} or Theorem 16.8 in
\cite{wo}). Granted this, and using also Corollary 4.3.3 in \cite{wo},
it is obvious that $(a)\sns (c)$. Now, $(b)\sns (c)$ according to
Corollary \ref{c2}.\qed
\begin{corol}
Let $A$ be a $\sigma$-unital (non-unital) purely infinite simple
\C. Then the units of $M(A)/A$ lift to units of $M(A)$ if and only if
$K_0(A)=0$.
\end{corol}
{\em Proof:}
By \cite{zh}, Theorem 1.2, $A$ is stable and has real rank zero, and
by Theorem 1.4 and Proposition 1.5 in \cite{cu}, $A$ has weak
cancellation. Thus the result follows from Corollary \ref{c3}.\qed
\espai
We close by remarking the fact that if $A$ is simple, $\sigma$-unital
(non-unital), with real rank zero and weak cancellation, then
$U(M(A))$ is connected (and hence $K_1(M(A))=0$). Indeed, if all
projections in $A$ are infinite, then $A$ is purely infinite simple,
hence stable (\cite{zh}, Theorem 1.2 (i)), and thus $U(M(A))$ is
connected. On the other hand, if there is a nonzero finite projection
$p\in A$, then $pAp$ has stable rank one, by Theorem 7.6 in
\cite{agop}, whence also $A$ has stable rank one, and then $U(M(A))$
is connected, by Lemma 3.3 in \cite{linscan}. Hence, for these
algebras, the units of $M(A)/A$ can be lifted to units of $M(A)$ if
and only if $U(M(A)/A)$ is connected.

\section*{Acknowledgements}

The author is grateful to Gert Pedersen for inviting him to the
Mathematics Institute of the University of Copenhagen, where this
work was carried out. It is also a pleasure to thank Ken Goodearl for
many helpful comments on an earlier draft of this paper.

\end{document}